\documentclass[pdflatex,sn-mathphys-num]{sn-jnl}

\usepackage{graphicx}%
\usepackage{multirow}%
\usepackage{amsmath,amssymb,amsfonts}%
\usepackage{amsthm}%
\usepackage{mathrsfs}%
\usepackage[title]{appendix}%
\usepackage[svgnames]{xcolor}%
\usepackage{textcomp}%
\usepackage{manyfoot}%
\usepackage{booktabs}%
\usepackage{algorithm}%
\usepackage{algorithmicx}%
\usepackage{algpseudocode}%
\usepackage{listings}%
\usepackage{float}

\usepackage{xstring}
\usepackage{comment}
\newcommand{\EatOneArg}[1]{}
\usepackage{array,booktabs,xurl}

\usepackage{amsmath,amssymb,array,longtable,booktabs}
\newcommand{\lf}[1]{{\ttfamily\small #1.lean}}
\newcommand{\br}{\allowbreak}

\theoremstyle{thmstyleone}%
\theoremstyle{thmstyletwo}%

\theoremstyle{thmstylethree}%

\begin{document}

\title[Certified neural computation of water waves]{Rigorously Certified
Neural-Network Computation of Nonlinear Water Waves at Near-Machine Precision}


\author*[1]{\fnm{Michal} \sur{Shavit}}\email{ms14479@nyu.edu}

\author[1]{\fnm{Yongji} \sur{Wang}}

\author[1]{\fnm{Tristan} \sur{Buckmaster}}

\affil*[1]{\orgdiv{Courant Institute of Mathematical Sciences}, \orgname{New
York University}, \city{NY}, \postcode{10012},  \country{USA}}




\abstract{
We present a computational framework that combines machine-learning-based discovery with rigorous computer-assisted numerical certification and formal verification for strongly nonlinear partial differential equations. We demonstrate the approach for periodic gravity–capillary traveling water waves. Physics-informed neural networks are used as a search mechanism, with the wave profile, wave speed and, in continuation experiments, surface tension learned directly from the governing equations without training data. Through a staged training and refinement procedure, the neural networks achieve near-machine-precision residuals of the governing equations, even in strongly nonlinear regimes. A single formulation recovers three classical high-amplitude regimes: pure-gravity Stokes waves, pure-capillary Crapper waves and gravity–capillary Wilton ripples. Starting from a high-amplitude Wilton ripple, continuation in surface tension further reveals new high-amplitude candidate solutions. We freeze one such machine-learned candidate as an explicit finite representation and certify its residual using interval arithmetic, rigorously accounting for floating-point roundoff, Fourier truncation, aliasing and grid-to-continuum errors. The resulting computation establishes a rigorous uniform bound on the governing-equation residual, independently validating the accuracy achieved through neural-network computation. We then formalize the mathematical structure of the certification argument in Lean, which machine-checks the deductions connecting the certified numerical bounds to the final residual theorem. Both the interval certification and formalization were developed with substantial assistance from a large language model, while the validity of the result remains independent of the language model: the numerical bounds are rigorously enclosed by interval arithmetic and the subsequent mathematical deductions are checked by the Lean kernel. Together, these steps provide a pipeline from neural discovery to rigorous numerical certification and machine-checked mathematical verification, enabling machine learning to explore strongly nonlinear regimes while retaining explicit and auditable mathematical guarantees.

}




\maketitle

\begin{center}
\href{https://github.com/michalshavitNYU/Rigorously-certified-water-waves}
{\texttt{Paper repository: computer-assisted accuracy proof and Lean verification}}
\end{center}

\tableofcontents

\section{Introduction}\label{sec1}

Machine learning can produce numerical solutions of nonlinear
partial differential equations with remarkably small
residuals at finitely many evaluation points~\cite{wang2024multi,wang2025sharpgradients}.
However, a small sampled residual alone does not provide a
rigorous bound on the continuous residual or on the error
relative to an exact
solution~\cite{mishra2023generalization,hillebrecht2022pde}.

Here we address continuous-residual certification by combining machine-learning-based discovery with rigorous numerical certification and formal verification. We demonstrate the approach for strongly nonlinear, high-amplitude water waves, where physics-informed neural networks discover candidate solutions with residuals approaching the limits of double-precision arithmetic. We then convert a selected candidate into a rigorous uniform residual certificate using interval arithmetic and formally verify the mathematical deductions underlying the certificate in Lean. The resulting pipeline provides a route from machine-learned numerical discovery to certified, machine-checkable scientific computation.

Surface gravity--capillary waves constitute a classical and enduring problem of fluid mechanics. Their governing equations are remarkably rich, describing phenomena ranging from small-amplitude capillary ripples to extreme ocean waves, yet are notoriously difficult to analyze because of their full nonlinearity and nonlocality. The subject has a long history spanning oceanography, mathematical physics and rigorous analysis, from Stokes' expansions~\cite{stokes1847theory} and Zakharov's Hamiltonian formulation~\cite{zakharov1968stability} to results on well-posedness and asymptotics~\cite{nalimov1974cauchy,lannes2013water}; see also the survey of steady water waves in Ref.~\cite{haziot2022traveling}. Despite this progress, the strongly nonlinear, high-amplitude regime remains substantially more difficult to access. In this regime, perturbative descriptions cease to apply, solution branches can develop intricate structure, and even the systematic construction and accurate numerical representation of nonlinear wave solutions become challenging.

This combination of strong nonlinearity, nonlocality and a rich solution structure makes water waves a stringent setting in which to test whether machine learning can move beyond numerical approximation and become a tool for mathematical discovery with rigorously certified performance. We address this question through three stages: neural discovery, computer-assisted certification by interval arithmetic, and formal verification.

First, we employ physics-informed neural networks (PINNs) \cite{raissi2019physics} as a search mechanism for high-amplitude two-dimensional periodic traveling solutions of the fully nonlinear water-wave equations. The network is optimized directly to satisfy the equations in conformal variables, with the wave profile, wave speed, and, in selected continuation experiments, the surface-tension parameter learned simultaneously. In this way, the network parametrizes and explores nonlinear solution branches. Building on the two-stage training procedure of ref.~\cite{wang2024multi,wang2025sharpgradients}, we obtain neural representations whose residuals approach the limits of double-precision arithmetic. Using the same formulation, we recover three canonical high-amplitude regimes: pure-gravity Stokes waves \cite{toland1996stokes}, pure-capillary Crapper waves \cite{crapper1957exact} and mixed gravity--capillary Wilton ripples \cite{wilton1915ripples}. Starting from a high-amplitude Wilton ripple and continuing in surface tension, the network further discovers new high-amplitude candidate solutions of the governing equation.

Second, we turn numerical accuracy into a rigorous statement using interval arithmetic. A small residual obtained by floating-point computation does not itself constitute a rigorous certificate, because roundoff, truncation, aliasing and discretization errors enter through the nonlinear and nonlocal operations in the water-wave equations. We therefore freeze a selected machine-learned candidate as an explicit finite representation and propagate rigorous intervals through the governing operators, explicitly accounting for floating-point roundoff, Fourier truncation, aliasing errors arising from grid-based evaluation, and the error incurred in passing from a finite grid to the continuous domain. This yields a certified uniform upper bound on the residual at near-machine precision, rather than a small residual observed only at finitely many numerical evaluation points.

Third, we formally verify the mathematical structure of the certification argument in Lean ~\cite{demoura2021lean}, a proof assistant and functional programming language, and the mathlib library~\cite{mathlib2020}. At this level of accuracy, the numerical certificate is itself necessarily computer assisted: our interval-arithmetic proof is implemented in Julia using the Arb library~\cite{johansson2017arb} and consists of hundreds of lines of code that propagate rigorous bounds through nonlinear and nonlocal operations. Although the resulting computation is mathematically rigorous, its implementation is too extensive and delicate to be audited line by line through conventional peer review. We therefore use Lean to provide a machine-checkable logical layer for the certification argument. Lean verifies the mathematical identities and deductions connecting the certified numerical bounds to the final residual theorem, substantially reducing the amount of computational reasoning that must be manually audited. The numerical interval bounds constitute the numerical trust boundary, while the subsequent mathematical deductions are checked by the Lean kernel.

Both the interval-arithmetic certification and its formalization were developed with substantial assistance from a large language model. Crucially, however, the validity of the resulting certificate does not depend on the reliability of the language model: the numerical bounds are rigorously established using interval arithmetic, while the mathematical deductions built upon them are checked by the Lean kernel. This separation allows a large language model to assist substantially in constructing a complex computer-assisted proof without becoming part of the trusted mathematical foundation of the result. The resulting code, analytic estimates, and formal deductions are documented through archived inputs, executable validation scripts, and Lean proofs, making the full certification pipeline available for inspection and reproduction.

Several lines of prior work motivate the individual components of this approach. PINN error analysis relates training residuals to solution errors under suitable stability assumptions~\cite{mishra2023generalization}, while a posteriori estimates have been developed for PDE-defined PINNs~\cite{hillebrecht2022pde}. Eiras et al.~\cite{eiras2024certification} developed methods for bounding PDE residuals over continuous domains rather than only at collocation points, and Tanaka and Yatabe~\cite{tanaka2026learnverify} combined PINN approximations with interval-based verification for nonlinear ordinary differential equations. High-accuracy neural training constitutes another ingredient, with natural-gradient and multistage methods reaching increasingly small residuals~\cite{muller2023highaccuracy,wang2024multi,wang2025sharpgradients}.

Rigorous numerics and proof assistants likewise have an established history. Computer-assisted analysis is surveyed in~\cite{gomezserrano2019computer}, while verified numerical computation has also been integrated with formal proof assistants~\cite{immler2018verified}. More recently, TorchLean~\cite{george2026torchlean} has developed formal infrastructure for neural-network semantics, differentiation and residual verification.

Together, the developments presented in this work establish a pipeline from machine-learning-based discovery to rigorous numerical certification and machine-checked mathematical verification. Pairs of these ingredients have previously been combined---including neural approximations with interval-based verification~\cite{tanaka2026learnverify} and rigorous numerics with formal verification~\cite{immler2018verified,george2026torchlean}---but, to our knowledge, the complete chain has not previously been closed from a solution discovered by a neural network, through rigorous interval certification, to formal verification of the resulting certificate in a proof assistant. In the present computation, floating-point roundoff, unresolved Fourier modes, aliasing and between-grid errors are enclosed within a single explicit bound for the continuous PDE residual. The framework therefore provides a route for using machine learning not only as a scientific-computing tool, but also as an instrument of mathematical discovery: exploring strongly nonlinear regimes that are difficult to access analytically while retaining explicit, rigorous and auditable mathematical guarantees.

\section{Results}\label{sec2}

Surface gravity-capillary waves are described in terms of the Euler equation for an incompressible, irrotational fluid of infinite depth with velocity potential $\phi$, velocity $\mathbf{v}=\nabla\phi$ and a free surface $y=\eta(x,t)$, where $x$ is the horizontal coordinate and $y$ is the vertical coordinate. Incompressibility gives the Laplace equation $\Delta\phi=0$ inside the fluid domain $\mathcal{D}=\{(x,y),\; -\infty<y\le \eta(x,t)\}$ with the boundary conditions at the surface $y=\eta(x,t)$:
\begin{align}\label{eq:ww}
  (\partial_t + \mathbf{v}\cdot\nabla)(y - \eta(x,t)) &= 0
 ,\\ 
  \phi_t+\tfrac12|\nabla\phi|^2+g\eta&=\kappa \frac{\eta_{xx}}{(1+\eta_x ^2)^{3/2}},
\end{align}
and $\nabla \phi \rightarrow 0$ as $y\rightarrow -\infty$. Here $g$ is the acceleration of gravity,  $\kappa$ is the surface tension coefficient and  
\begin{equation}
    \frac{\eta_{xx}}{(1+\eta_x ^2)^{3/2}}=\partial_x(\frac{\eta_{x}}{\sqrt{1+\eta_x ^2}})
\end{equation}
is the mean curvature of the free surface. 
The first boundary condition is the kinematic condition, expressing that the free surface moves with the fluid particles. The second is Bernoulli's condition, in which the pressure jump across the free surface is proportional to the mean curvature.

\begin{figure}[H]
\centering
\includegraphics[width=0.7\textwidth]{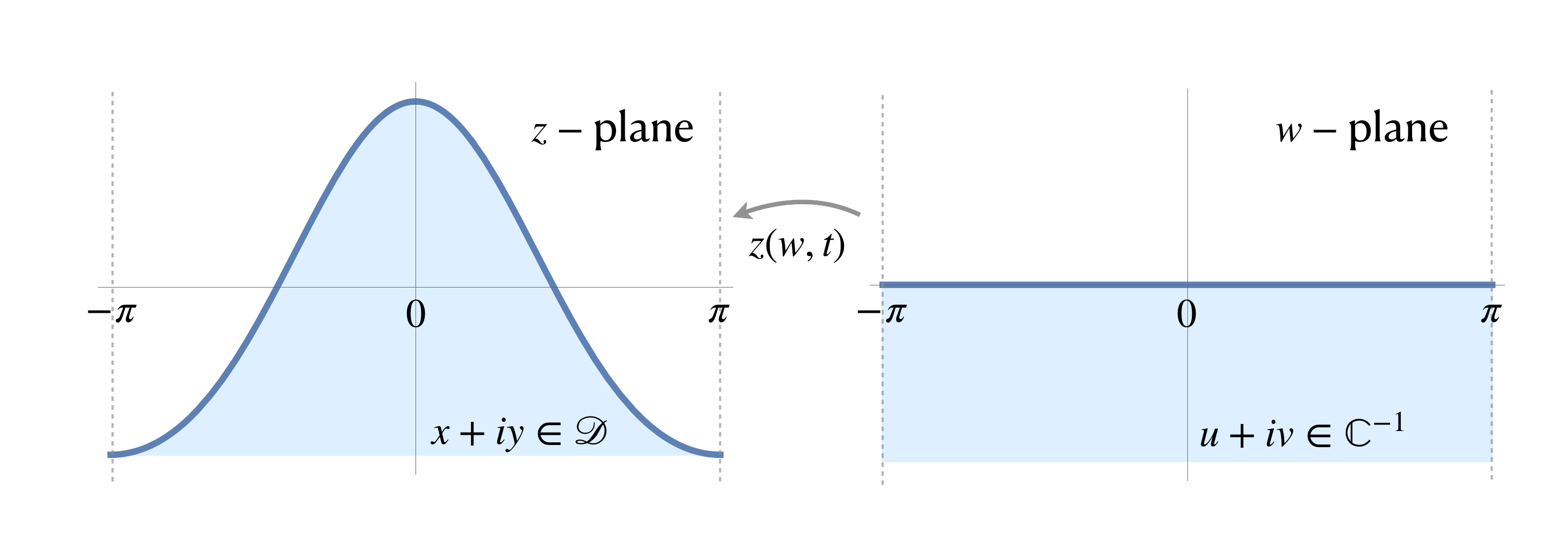}
\caption{\label{fig:complex} The conformal mapping from the region in the $w$ plane $(u,v)\in [-\pi,\pi]\times (-\infty,0]$ to the domain occupied by the fluid in the $(x,y)$ plane, with $(x,y)\in [-\pi,\pi]\times(-\infty,\eta(x,t)]$. The line $v=0$ (blue line) is mapped onto the free surface $\eta(x,t)$.}
\end{figure}

For periodic waves traveling in one direction, one can reformulate the equations above in terms of a conformal map $z(w,t):=x(w,t)+iy(w,t)$ that maps one periodic strip of the lower half-plane $\{(u,v)\in[-\pi,\pi]\times(-\infty,0]\}$ to the fluid domain $\mathcal{D}$ as illustrated in Fig.~\ref{fig:complex}. 

Conformal formulations are described in ~\cite{tanveer1991singularities,dyachenko1996nonlinear,haziot2022traveling}. On $v=0$, steady traveling waves $y=y(u-ct)$ satisfy a gravity--capillary extension of the Babenko formulation~\cite{babenko1987remarks,haziot2022traveling}:

\begin{align}\label{eq:babenko}
  \mathcal{B}(y) =-\mu  \hat H y_u+ & \left(y + y\,\hat H y_u + \hat H(yy_u)\right) \\
   - &\beta\,\partial_u\!\left(\frac{y_u}{\sqrt{y_u^2+(1+\hat H y_u)^2}}\right)
   + \beta\,\hat H\partial_u\!\left(\frac{1+\hat H y_u}{\sqrt{y_u^2+(1+\hat H y_u)^2}}\right)&=0 \nonumber,
\end{align}
where $\mu$ is the squared (rescaled) wave speed, $\beta=\kappa /(g \lambda^2)$ is the Bond number encoding surface tension, gravity, and wavelength $\lambda$. $y(u)$ is the elevation of the free surface in conformal coordinates and $\hat H$ is the periodic Hilbert transform defined on $L_2[-\pi,\pi]$ by linearity from the relation
\begin{equation}
    \hat{H}e^{inu}=-i\text{sign}(n)e^{inu}.
\end{equation}

The underlying water-wave equations have classical high-amplitude families in the limits of zero capillarity and zero gravity. The latter limit requires a separate scaling, since the positive-gravity normalization in Eq.~\eqref{eq:babenko} divides by $g$. At $\beta=0$, the classical one-parameter family of Stokes waves is represented as a curve in the speed and amplitude plane $(c,a)$ connecting small amplitude solutions with a limiting Stokes wave of extreme height with a $2\pi/3$ crest angle \cite{toland1996stokes}. The amplitude is defined as the difference between the maximal and minimal height divided by the period
\begin{equation}
    a = \frac{\text{max}y-\text{min}y}{2\pi}.
\end{equation}
At zero gravity, $g=0$, there is a family of pure capillary solutions, found by Crapper~\cite{crapper1957exact}, which is represented as a curve in the speed and amplitude plane $(c,a)$ connecting small amplitude capillary waves with a limiting Crapper wave of extreme height that encloses a bubble. 

When both gravity and capillarity are present, resonant families of solutions, known as Wilton ripples, arise at special values of the surface tension. At the linear level, a Wilton ripple contains two harmonics, \(\cos(mu)\) and \(\cos(nu)\), whose phase speeds coincide. While small-amplitude Wilton ripples are well understood and have been computed numerically, their strongly nonlinear, high-amplitude regime remains much less explored, both analytically and numerically. More broadly, numerical studies have long investigated finite- and large-amplitude gravity--capillary waves~\cite{schwartz1979capillary,chen1980capillary}, and modern work has examined the spectral stability of small-amplitude Wilton ripples~\cite{trichtchenko2016wilton}. The structure and continuation of high-amplitude Wilton branches remain challenging; even the choice of amplitude constraint can affect access to connected families~\cite{shelton2021parasitic}. More generally, the gravity--capillary problem has a rich bifurcation structure in the wave-speed--surface-tension parameter space. Local and secondary bifurcations, including resonant Wilton branches associated with multiple eigenvalues of the linearized problem, have been studied analytically and numerically~\cite{toland1985bifurcation,jones1986symmetry,buffoni2000regularity}, with Babenko-type formulations providing a useful framework for such analysis~\cite{buffoni2000regularity,wahlen2023global}. Yet the global organization of strongly nonlinear, high-amplitude branches remains incompletely characterized. In this setting, a rigorous residual bound is already of independent value as a quantitative diagnostic, even without the additional existence argument required to prove a nearby exact solution. 

\subsection{Neural discovery and recovery of high-amplitude water waves}

\begin{figure}[h!]
\centering
\includegraphics[width=\textwidth]{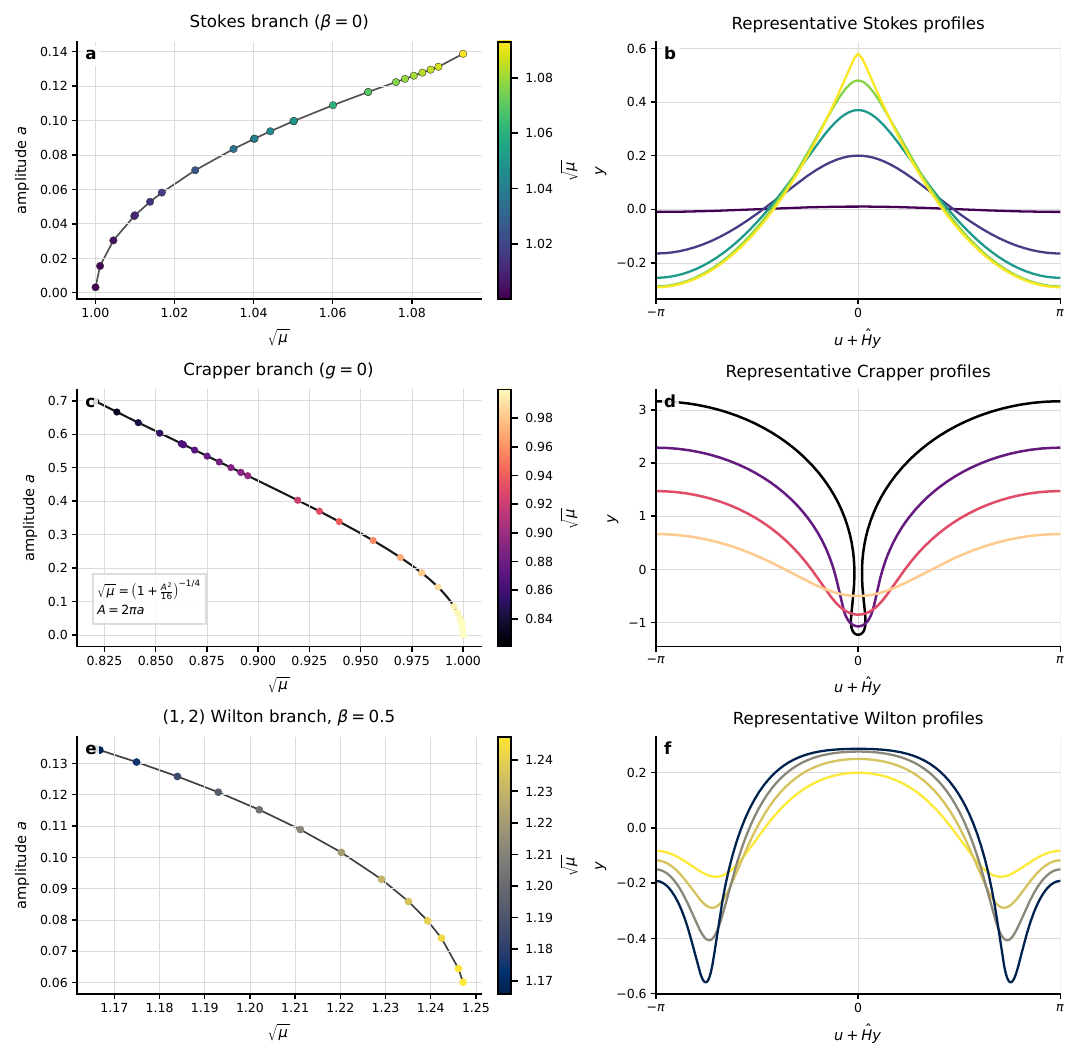}
\caption{\label{fig:benchmark_regimes}\textbf{Three high-amplitude regimes reproduced by the PINN.} 
\textbf{a,} Pure-gravity Stokes branch with $\beta=0$, plotted in the $(\sqrt{\mu},a)$ plane, where $a=(\max y-\min y)/(2\pi)$. The panel shows the selected Stokes files, including the low-amplitude main-branch runs, that satisfy $\max|\mathcal{B}(y_N)|\leq 10^{-7}$, with the grey curve joining the retained runs. 
\textbf{b,} Five representative retained Stokes profiles plotted over $[-\pi,\pi]$ as $x(u)=u+\hat H y(u)$ versus $y(u)$. 
\textbf{c,} Pure-capillary Crapper branch with $g=0$. All Crapper files in the plotting folder are included as points on the branch, and the black curve shows the displayed relation $\sqrt{\mu}=(1+A^2/16)^{-1/4}$ with $A=2\pi a$. 
\textbf{d,} Representative Crapper profiles plotted over $[-\pi,\pi]$. 
\textbf{e,} Reproduced $(1,2)$ Wilton branch at $\beta=0.5$, again with all available branch files included. 
\textbf{f,} Four representative Wilton free-surface profiles, including low- and high-amplitude cases, plotted over $[-\pi,\pi]$.}
\end{figure}

We use physics-informed neural networks (PINNs) \cite{raissi2019physics} to discover water-wave profiles and their associated physical parameters directly from the Babenko equation \eqref{eq:babenko}. Rather than being trained on precomputed or experimental data, the network is optimized against the governing equation together with the physical and normalization constraints defining the solution family under consideration. The same neural architecture and residual-based optimization framework is used throughout, while the physical parameters and auxiliary constraints are adapted to select the relevant branch. The trained profile is denoted by $y_P$. For each run, we minimize a weighted squared residual of the Babenko equation and, when used, its differentiated form, together with auxiliary conditions controlling quantities such as vertical translation, amplitude and symmetry. The wave speed $c_P=\sqrt{\mu_P}$ is optimized simultaneously with the profile. Depending on the experiment, gravity and surface tension are fixed to select a particular physical regime, while in the continuation experiments the surface-tension parameter is itself allowed to vary. Details of the network architecture, Hilbert-transform evaluation, loss function, constraints and optimization procedure are given in Methods.

We first use three canonical high-amplitude solution families as benchmarks for the discovery framework (Fig.~\ref{fig:benchmark_regimes}). These regimes probe qualitatively different balances in the same nonlinear equation. In the zero-surface-tension limit, $\beta=0$, the network reproduces the classical pure-gravity Stokes branch \cite{stokes1847theory}. In the zero-gravity limit, $g=0$, it recovers the pure-capillary Crapper family \cite{crapper1957exact}. When both gravity and surface tension are present, the network recovers high-amplitude $(1,2)$ Wilton ripples \cite{wilton1915ripples}, for which the first and second Fourier harmonics travel at the same speed at the linear level. The three families have markedly different geometries, dominant physical balances and Fourier structures, and their recovery therefore provides complementary tests of the same computational framework. The physical and normalization constraints differ between these regimes, while the neural representation and equation-based optimization remain common.

The Wilton regime provides a particularly demanding test. At small amplitude, the $(1,2)$ Wilton branch originates from a linear resonance between two Fourier modes: $\cos(u)$ and $\cos(2u)$, but at finite amplitude, their contributions are coupled by the full nonlinear water-wave equation. The amplitudes of the resonant components therefore cannot in general be varied independently: nonlinear constraints determine how the harmonics combine along the solution branch. Continuing a high-amplitude Wilton ripple is consequently not equivalent to simply varying the amplitude 
parameter and reoptimizing an unconstrained profile. The continuation must simultaneously preserve the appropriate nonlinear constraints while allowing the profile and wave speed parameter to adjust, making the high-amplitude resonant regime a substantially more restrictive search problem.


\subsubsection{Continuation experiments and new solutions}
The benchmark and continuation experiments separate two roles of the neural computation. Recovery of the Stokes, Crapper and Wilton families tests whether the framework can represent and optimize solutions across distinct nonlinear regimes when supplied with the appropriate physical constraints. Continuation away from the high-amplitude Wilton state then uses the same machinery in a discovery setting, where the resulting profiles are not prescribed in advance. The low-residual candidates obtained in this stage provide the starting point for the rigorous certification described below.

\begin{figure}[h!]
\centering
\includegraphics[width=\textwidth]{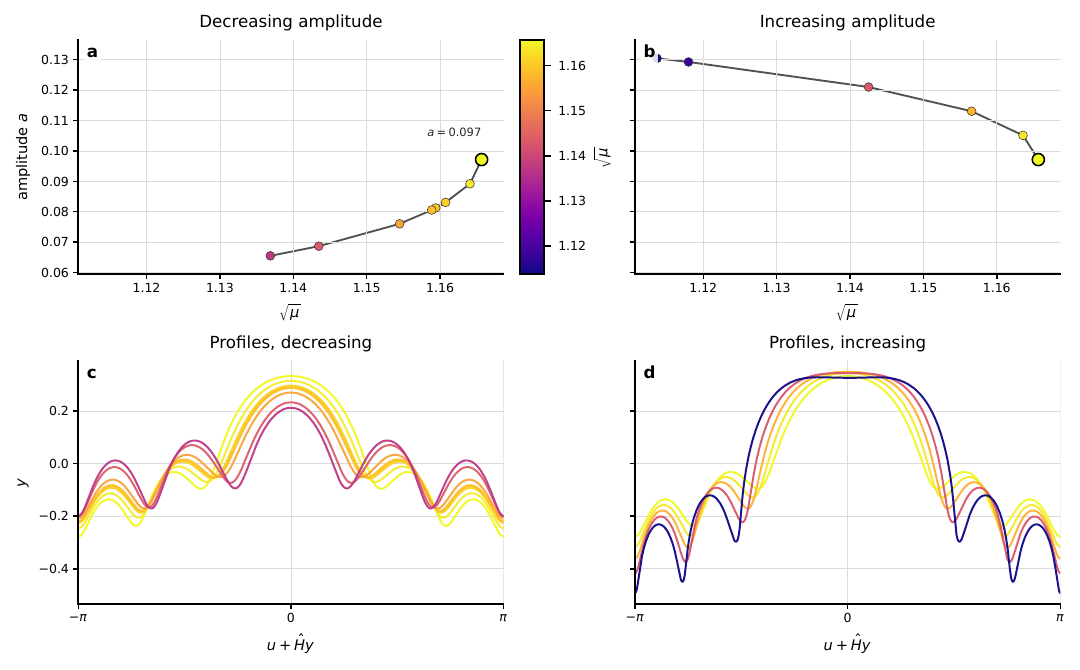}
\caption{\label{fig:wilton_continuation_s}\textbf{Continuation away from the high-amplitude $(1,2)$ Wilton ripple.}
\textbf{a,} Amplitude-decreasing continuation at $\beta=0.498^2\sim 0.248004$, starting from the shared seed with $0.09<a<0.1$ and plotted in the $(\sqrt\mu,a)$ plane. 
\textbf{b,} Amplitude-increasing continuation from the same seed point. 
\textbf{c,} Free-surface profiles along the amplitude-decreasing branch. 
\textbf{d,} Free-surface profiles along the amplitude-increasing branch, with colours in both profile panels matched to the same $\sqrt{\mu}$ scale.}
\end{figure}

After reproducing a high-amplitude $(1,2)$ Wilton ripple, we use it as a seed for neural continuation. We allow the surface-tension parameter to vary together with the wave profile and speed. Starting from the same high-amplitude seed, continuation in different constrained directions produces new gravity--capillary profiles with small Babenko residuals (Figs.~\ref{fig:wilton_continuation_s} and~\ref{fig:wilton-beta-continuation-four}). Although these solutions are reached by continuation from a Wilton ripple, they are not simply nearby members of the known high-amplitude Wilton family: their Fourier composition and values of $\beta$ depart substantially from those of the corresponding Wilton branch, and their geometry develops qualitatively different features as the continuation proceeds. Toward one end of the computed branch, the crest develops a pronounced central depression while the trough becomes increasingly sharp, producing profiles distinct from the high-amplitude Wilton ripples used to initialize the search. To our knowledge, gravity--capillary solutions with this combination of Fourier structure, parameter values and free-surface geometry have not previously been computed. For comparison with the established large-amplitude
gravity–capillary families, see ~\cite{schwartz1979capillary,chen1980capillary,shelton2021parasitic}.

\begin{figure}[h!]
    \centering
    \includegraphics[width=\linewidth]{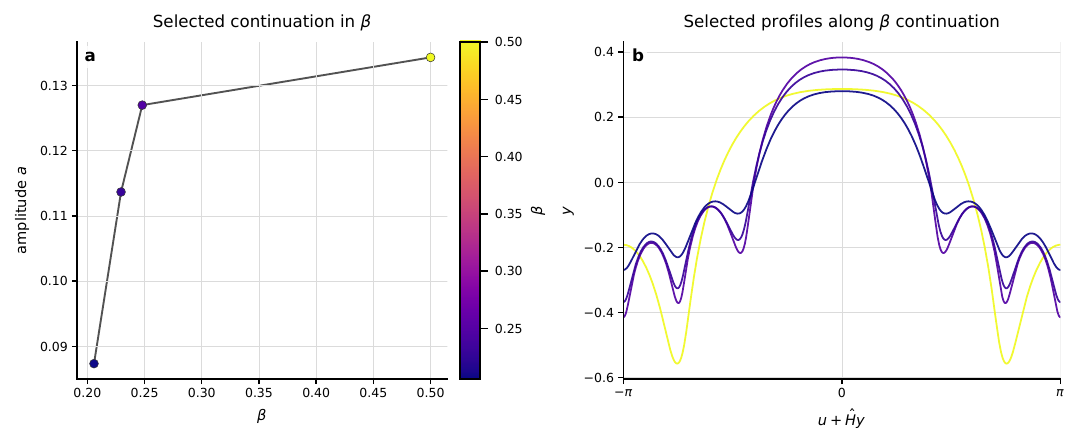}
    \caption{
    Continuation near a high-amplitude Wilton ripple.
    \textbf{a}, Selected points on the continuation branch, shown in the amplitude--$\beta$ plane.
    \textbf{b}, Corresponding free-surface profiles over one period, plotted against the physical horizontal coordinate $x(u)=u+\hat H y(u)$ and colored by the same $\beta$ values as in \textbf{a}.
    The separately displayed degree-150 polynomial in Fig.~\ref{fig:wilton-certified-profile} is the candidate used for the final residual certificate.
    }
    \label{fig:wilton-beta-continuation-four}
\end{figure}

\newpage

We were also unable to continue this family toward the small-amplitude regime. Within the range accessible to our numerical continuation, the branch instead remains intrinsically finite amplitude. This observation suggests that these profiles may not arise through a conventional bifurcation from a linear water-wave mode, although establishing whether the family is genuinely disconnected from the small-amplitude spectrum requires further analysis. In this sense, the high-amplitude Wilton ripple provides an entry point into a region of the nonlinear solution landscape that is not apparent from the linearized problem.

A detailed mathematical characterization of this solution family, including its global bifurcation structure, possible limiting states, stability and connection, or lack thereof, to small-amplitude solutions, is beyond the scope of the present work. Here these solutions primarily demonstrate the exploratory capability of the computational framework: constrained neural continuation can identify previously uncharacterized high-amplitude states, which can subsequently be passed to rigorous numerical certification. Our focus in the remainder of the work is therefore on establishing this discovery-to-certification pipeline and on rigorously validating a representative machine-learned solution.

\subsection{Second-stage training and a near-machine-precision candidate}

We demonstrate that the Babenko residual can be reduced to near
the double-precision round-off level. Successive neural corrections,
followed by Fourier projection and Newton polishing, reduced the
maximum sampled scalar residual by approximately seven orders of
magnitude, from $\mathcal{O}(10^{-7})$ to $\mathcal{O}(10^{-14})$, see Fig.~\ref{fig:residual-refinement}. Independent interval validation
then established a rigorous uniform residual bound of
$\mathcal{O}(10^{-13})$ for the final polynomial.

\begin{figure}[h!]
    \centering
    \includegraphics[width=\linewidth]{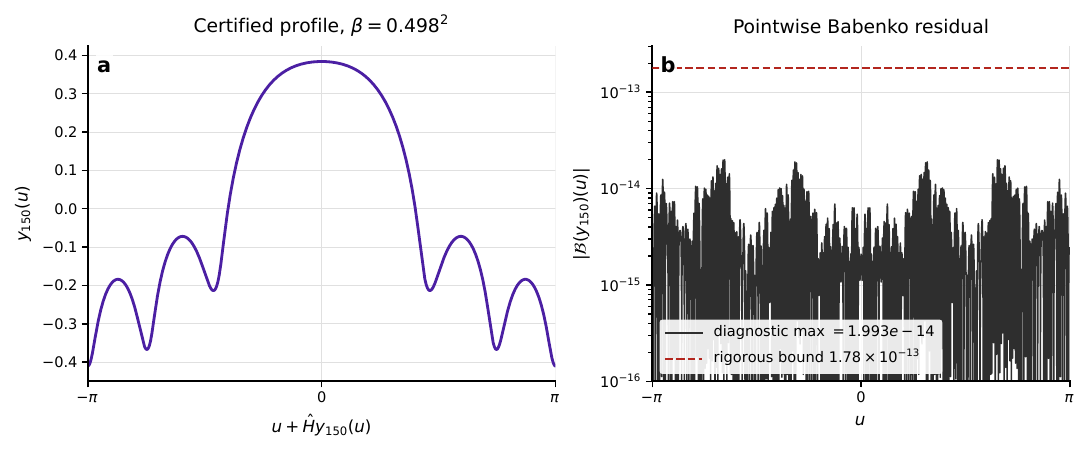}
    \caption{
    Certified high-amplitude Wilton-ripple profile.
    \textbf{a}, Free-surface profile for the certified solution at $\beta=0.498^2$, plotted over one period against the physical horizontal coordinate $u+\hat{H}y$.
    \textbf{b}, Pointwise residual $|\mathcal{B}(y_N)|$ for the same solution. The dashed line marks the certificate level $\|\mathcal{B}(y_N)\|_{\infty} \leq 1.78\times10^{-13}$ with $N=150$.
    }
    \label{fig:wilton-certified-profile}
\end{figure}

\subsection{Certified residual estimates and formal verification in Lean}

We now turn from numerical discovery to rigorous certification. The
neural-network approximation is used to construct a candidate profile;
the object ultimately certified is an explicit finite trigonometric
polynomial, with fixed coefficients and parameters.

Let
\begin{equation}
y_N(u)=\frac{a_0}{2}+\sum_{k=1}^{N}a_k\cos(ku),
\end{equation}
and let \(\mu,\beta>0\) be fixed rational parameters. We write
\begin{equation}
\mathcal{R}_N(u)=\mathcal{B}_{\mu,\beta}(y_N)(u)
\end{equation}
for the residual of the gravity--capillary Babenko equation. We establish the rigorous uniform estimate
\begin{equation}
\|\mathcal{R}_N\|_{L^\infty(\mathbb{R})}
=
\sup_{u\in\mathbb{R}}|\mathcal{R}_N(u)|
\leq \varepsilon,
\label{eq:certification-goal}
\end{equation}
for an explicitly calculated constant \(\varepsilon\). For the candidate certified
here,
\begin{equation}
N=150,
\qquad
\varepsilon=\frac{178}{10^{15}}
=1.78\times10^{-13}.
\label{eq:certified-tolerance}
\end{equation}
This is a uniform residual theorem for a supplied finite series. It is not an existence or uniqueness theorem for an exact nearby water wave. Such a result would additionally require a rigorous inverse bound for the appropriately constrained linearized Babenko operator and control of the nonlinear remainder \cite{gomezserrano2019computer}; these ingredients are not part of the present certificate.

The certification has two complementary components: validated numerical computation and formal verification. Julia/Arb performs the substantive interval-arithmetic calculation, rigorously enclosing numerical errors and producing certified bounds. Lean~4 formalizes the analytic framework and verifies the deductions from these certified numerical inputs to the final global theorem.

This certification has both practical numerical and pure-mathematical motivations. From a numerical perspective, the objective is to replace an empirically observed small residual by a rigorous uniform bound at the same scale. Numerically observing a small residual and rigorously certifying such a bound are distinct achievements.

From a pure-mathematical perspective, the objective is to make
PINN-based discovery usable as a high-precision tool for studying
highly nonlinear and nonlocal equations. The resulting computations
can serve as quantitative inputs to further mathematical analysis:
the candidate is explicitly defined, its residual is rigorously
controlled, and the computational and logical trust boundaries are
made explicit. In particular, the residual certificate provides
a first rigorous ingredient toward a computer-assisted proof of the
existence of a nearby exact solution.

Source code, archived run evidence, build instructions, and the supplementary mathematical account are indexed in the
\href{https://github.com/michalshavitNYU/Rigorously-certified-water-waves}{paper repository}.

We now describe the two stages in turn.

\subsubsection{Rigorous certification by interval arithmetic}

The interval calculation is implemented in the Julia package, which produces a computer-assisted proof
certificate bounding the residual of Babenko's equation for steady
periodic water waves. Given a candidate profile represented as a
truncated even cosine series, the package rigorously evaluates the
Babenko residual and generates a Lean~4 certificate encoding the certified numerical bounds as exact rational literals, ready for verification by Lean.

The interval calculation rigorously accounts for floating-point
roundoff, Fourier truncation and aliasing in nonlinear terms, and the
passage from grid-based estimates to bounds valid throughout the
continuous domain. For the finite trigonometric profile \(y_N\) and
its derivatives, the periodic Hilbert transform is evaluated exactly
mode by mode. For nonlinear quantities in the capillary term, which
are not finite trigonometric polynomials, unresolved Fourier modes and
aliasing errors are bounded explicitly.

The computation also establishes a strictly positive lower bound for
\[
D(u)=
\sqrt{
\bigl(\partial_u y_N(u)\bigr)^2+
\bigl(1+\hat H\partial_u y_N(u)\bigr)^2
},
\]
ensuring that the capillary denominator remains nonzero. Rigorous
Taylor bounds on each grid cell then promote the discrete interval
estimates to the full period and, by periodicity, to all
\(u\in\mathbb{R}\). Thus,
\begin{equation}
\|\mathcal{B}_{\mu,\beta}(y_N)\|_{L^\infty(\mathbb{R})}
\leq \varepsilon,
\qquad
D(u)\geq d_{\min}>0 .
\end{equation}

The implementation, input coefficients, run logs and reproduction
commands are available in the
\href{https://github.com/michalshavitNYU/Rigorously-certified-water-waves/tree/main/Julia}
{Julia directory}.

\subsubsection{Formal verification in Lean}

The mathematical verification layer of the computer-assisted proof is implemented in Lean~4. Julia/Arblib validates the fixed cosine-series profile, rigorously controlling rounding, Fourier truncation, aliasing and between-grid errors, and exports the resulting certified bounds together with the profile and computation parameters as exact rational data. A reusable Lean library formalizes the mathematical structure underlying the certification, including the periodic Hilbert transform, the finite trigonometric profile, the capillary terms and the Babenko residual. It proves the mode-by-mode action of the Hilbert transform on finite trigonometric polynomials, the Fourier-tail and aliasing estimates required for the nonlinear capillary terms, and the derivative, Taylor-remainder and grid-to-continuum estimates used in the certification.

For the candidate considered here, Lean takes the certified denominator bound and one-period residual bound produced by Julia/Arblib as explicit inputs, verifies the analytic deductions and periodic extension to the full real line, and checks the exact rational comparisons yielding the reported constants. Thus, validated numerics establishes the numerical inequalities, while Lean verifies the mathematical deductions from them. The resulting formal library is independent of the particular numerical profile and can be reused for future certified solutions.

The
\href{https://github.com/michalshavitNYU/Rigorously-certified-water-waves/tree/main/Lean}
{Lean directory}
contains the formal definitions, theorem modules,
certificate deductions, axiom audit and pinned build instructions.

\section{Discussion}\label{sec:discussion}

The results presented here demonstrate a route from machine-learning-based exploration of a nonlinear PDE to a rigorously certified and machine-checkable computational statement. The role of the neural network is primarily one of discovery: it provides a flexible representation with which to search nonlinear solution spaces and identify high-accuracy candidate solutions. The subsequent certification is deliberately separated from the training procedure. Once a candidate is selected, it is frozen as an explicit finite representation, and the mathematical claim concerns this representation rather than the neural network, optimizer or training history. This separation allows the exploratory flexibility of machine learning to be combined with the reproducibility and explicit error control of validated numerical computation.

The water-wave problem illustrates both the potential and the limitations of this approach. The same computational framework accesses qualitatively different high-amplitude regimes, including pure-gravity, pure-capillary and resonant gravity--capillary waves, but the search is not free of problem-specific mathematical information. Different solution families require appropriate physical, normalization and symmetry constraints, and continuation of high-amplitude Wilton ripples requires additional control of the relative amplitudes of the resonant Fourier modes. The neural network therefore does not replace mathematical structure; rather, it provides a flexible computational representation within which such structure can be imposed while leaving the detailed nonlinear profile to be determined by the governing equation. Developing methods that can infer useful continuation coordinates and constraints directly from the learned solution landscape would reduce this dependence on prior analytical knowledge.

A second distinction concerns numerical accuracy and mathematical existence. The certificate established here is a rigorous uniform bound on the residual of an explicit approximate solution. It shows that the reported near-machine-precision residual is not an artifact of floating-point roundoff, Fourier truncation, aliasing or finite sampling. It does not, by itself, prove the existence of an exact water-wave solution nearby. Such a result would require quantitative control of the linearized Babenko operator and of the nonlinear remainder, for example through a Newton--Kantorovich or related contraction argument. The certified residual obtained here provides one of the quantitative inputs required for such an existence proof and establishes a natural interface between machine-learning-based discovery and computer-assisted analysis.

Formal verification addresses a different source of uncertainty. At the accuracy considered here, the interval-arithmetic certification is itself a substantial computation involving hundreds of lines of validated numerical code and a sequence of analytic estimates that would be difficult to audit completely by conventional inspection. Lean makes the mathematical identities and deductions connecting the numerical certificate to the final theorem machine-checkable. The present implementation nevertheless retains an explicit numerical trust boundary: Lean does not independently reproduce the Arb computation, and the validated numerical bounds supplied by the interval-arithmetic calculation remain trusted inputs to the formal layer. 

An important feature of the workflow is that substantial parts of both the rigorous numerical implementation and the formalization were developed with assistance from a large language model. This illustrates a potentially broader role for generative AI in scientific computing: not only in producing numerical candidates, but also in helping construct the machinery required to validate them. The reliability of the final mathematical statement need not inherit the reliability of the generative model, provided that model-generated computations and arguments terminate in independently checkable certificates. In the present framework, interval arithmetic provides rigorous numerical enclosures and Lean checks the subsequent logical deductions.

More broadly, the approach is not specific to water waves or to PINNs. Its essential ingredients are a computational method capable of discovering accurate candidate solutions, an explicit representation to which validated numerics can be applied, rigorous control of discretization and continuum errors, and a formal layer capable of checking the deductions made from the resulting numerical bounds. These ingredients suggest a general architecture for trustworthy machine-learning-assisted scientific computing, particularly for nonlinear problems in which analytical descriptions are incomplete and numerical exploration plays an important role. 

\section{Methods}\label{sec:methods}

\subsection{PINN architecture}

\subsubsection{Neural representation}

The unknown free surface is represented in conformal coordinates by a neural function of the periodic variable $u\in[-\pi,\pi]$. For the symmetric runs, periodicity and evenness are imposed by feeding the network the transformed coordinate $\cos u$. The implementation used for the plotted runs consists of three scalar
fully connected multilayer perceptrons with tanh activations,
initialized with Xavier-scaled truncated normal
weights~\cite{glorot2010understanding}. Each scalar network has three hidden layers of width 30. The networks are combined as
\begin{equation}
  y_P(u)=s\,N_0(\cos u),\qquad
  T_{1,P}(u)=\,N_1(\cos u)\sin u,\qquad
  T_{2,P}(u)=\,N_2(\cos u),
\end{equation}
where $s$ is the amplitude scale used in a given continuation step. This construction makes $y_P$ and $T_{2,P}$ even and $T_{1,P}$ odd, matching the crest-symmetric geometry. The auxiliary variables $T_{1,P}$ and $T_{2,P}$ represent the two components of the normalized tangent/curvature expression,
\begin{equation}
  T_1=\frac{y_u}{\sqrt{y_u^2+(1+\hat H y_u)^2}},
  \qquad
  T_2=\frac{1+\hat H y_u}{\sqrt{y_u^2+(1+\hat H y_u)^2}},
\end{equation}
so that the surface-tension term can be evaluated stably through first derivatives of the auxiliary outputs.

The wave speed is optimized together with the network weights through a bounded scalar parameterization
\begin{equation}
  c_P=c_0+\eta_c\tanh\rho_c,
\end{equation}
where $c_0$ is the continuation seed and $\eta_c$ is the allowed local search radius. In the Wilton continuation runs in which surface tension is trained, the same bounded parametrization is used for $\beta_P$. Fixed-regime benchmark runs set the corresponding physical parameter to its prescribed value.

\subsubsection{Residual, Hilbert transform and loss}

All derivatives of the neural outputs are obtained by automatic
differentiation in JAX~\cite{jax2018github}, with double-precision
arithmetic enabled. The periodic Hilbert transform is evaluated numerically on a dense nonuniform grid of $15,000$ sample points. The singular integral is split into a local contribution around the evaluation point, a close-field correction and an outer contribution evaluated by Gaussian quadrature. The same Hilbert-transform routine is used for $y_P$, $y_{P,u}$, derivatives of $y_P^2$ and derivatives of the auxiliary variables entering the surface-tension term.

Training minimizes a weighted loss
\begin{equation}
  \mathcal{L}
  =
  w_{\rm eq}\sum_{u_i\in\mathcal{C}}\omega_i
  \left\| \mathcal{R}_P(u_i)\right\|^2
  + w_{\rm der}\sum_{v_j\in\mathcal{D}}\nu_j
  \left|\partial_u\mathcal{B}(y_P)(v_j)\right|^2
  + \mathcal{L}_{\rm aux}.
\end{equation}
Here $\mathcal{R}_P=(e_1,e_2,e_3)$, where $e_1$ is the main Babenko residual and \(e_2\) and \(e_3\) enforce the definitions of \(T_{2,P}\) and \(T_{1,P}\), respectively. The auxiliary loss $\mathcal{L}_{\rm aux}$ fixes the amplitude and, when needed, removes vertical translation through a mean-zero condition or enforces or breaks prescribed symmetry conditions. The active runs use $w_{\rm eq}=0.1$ with residual-component weight one; the normalization and mean constraints, when used, are kept on throughout training.

\subsubsection{Collocation and optimization}

The collocation sets are resampled during optimization. In the runs used here, each optimization step combines fixed and moving collocation points: 251 fixed equation points, 201 fixed derivative points, two groups of 71 points concentrated near the endpoints and 121 moving points for each of the equation and derivative residuals. The sampling density is biased toward regions where high-amplitude waves develop sharp features, especially near crests and troughs. After a training pass, the point distribution can be updated from the residual itself so that later passes concentrate effort where the current approximation is least accurate.

We use the Kronecker-factored approximate curvature
(K-FAC) optimizer~\cite{martens2015kfac}, implemented in JAX,
with adaptive learning rate, momentum and damping,
a minimum damping of $10^{-14}$, and a typical training horizon
of $20{,}000$ iterations per continuation step. New branches are followed by warm starting from a previously trained profile and allowing the speed, and when appropriate the surface tension, to move within a prescribed local window. 

\subsection{Second-stage training}
\label{sec:second-stage-training}

Following the multistage refinement approach
of Refs.~\cite{wang2024multi,wang2025sharpgradients},
we refine the first-stage candidate through successive additive
neural corrections targeting its remaining residual. Let \(Y_p=(y_p,T_{1,p},T_{2,p})\) denote the
previously trained fields, frozen at the start of a refinement pass.
The corrected approximation is
\begin{equation}
Y_P(u)=Y_p(u)+\varepsilon_p\Phi_\theta(u),
\end{equation}
where \(Y_P=(y_P,T_{1,P},T_{2,P})\), \(\Phi_\theta\) is a
three-component correction network, and \(\varepsilon_p\) sets its
amplitude scale. Write \(\mathcal R(Y)\) for the three-field Babenko
residual operator and
\(\mathcal R_P=\mathcal R(Y_P)=(e_1,e_2,e_3)\). When those
definitions hold, \(e_1=\mathcal B(y_P)\), where \(\mathcal B\) is
the scalar Babenko operator in Eq.~\eqref{eq:babenko}.

The correction ansatz is motivated by linearization about the frozen
fields \cite{wang2024multi,wang2025sharpgradients}:
\begin{equation}
\mathcal R(Y_p+\varepsilon_p\Phi_\theta)
\simeq \mathcal R(Y_p)
+\varepsilon_p D\mathcal R(Y_p)\Phi_\theta.
\end{equation}
To cancel the existing residual to first order, the correction should
therefore approximately satisfy
\begin{equation}
\varepsilon_p D\mathcal R(Y_p)\Phi_\theta
\simeq -\mathcal R(Y_p).
\end{equation}
The linearized equation explains the correction strategy;
optimization evaluates the full nonlinear residual \(\mathcal R_P\).

The correction networks used sinusoidal activations and a tunable
first-layer frequency scale. We chose the frequency scale of the
correction to \(y_P\) from the prominent frequencies in the
remaining residual, and chose \(\varepsilon_p\) according to its
magnitude. These choices guide the scales represented by the network;
they do not constrain its output to a single Fourier mode or make
\(\varepsilon_p\) equal to a residual norm. The input \(\cos u\)
makes the outputs periodic and even, while a factor of \(\sin u\)
imposes odd parity on the \(T_{1,P}\) correction. After one correction
was trained, its contribution was incorporated into the frozen base
and a further, smaller-scale correction was trained. The capillary
parameter remained fixed, while the wave speed was allowed a small
adjustment.

The refinement stage evaluated the periodic Hilbert transform by FFT
on full uniform grids rather than by the nonuniform quadrature used
during discovery. Neural-field derivatives were computed by automatic
differentiation. Shifting the uniform grid by fractions of a grid
spacing sampled the residual between earlier collocation points
while preserving the periodic grid required by the FFT. Later
refinement passes used finer grids.

\subsubsection{Gradient-normalized residual weighting}

The second-stage refinement also employs gradient-normalized
residual weighting~\cite{wang2025sharpgradients}, adapted here
to the three-field Babenko system. A standard mean-square residual,
\[
\sum_j |\mathcal R(Y_P)(u_j)|^2,
\]
weights the residual in absolute terms. When the solution contains
steep or rapidly varying regions, however, the natural scale of the
equation can vary strongly across the domain, causing such regions to
dominate the optimization. We therefore normalize the residual by a
local scale determined from spatial derivatives of the frozen base
fields, leading schematically to
\[
\sum_j
\left|
\frac{\mathcal R(Y_P)(u_j)}
{\text{local derivative scale at }u_j}
\right|^2.
\]
In this sense, the loss measures a locally scaled residual rather than
a purely absolute one.

For the Babenko system, the second-stage loss retained the
component-wise weighted residual structure used during training.
Let \(\mathcal C\) denote the collocation points for \(\mathcal R_P\)
and \(\mathcal D\) those used for the derivative of the main equation.
To balance regions with different levels of spatial variation, the
residuals were divided by scales \(S_{\rm eq}\) and \(S_{\rm der}\)
constructed from spatial derivatives of the frozen fields \(Y_p\).
These scales combine derivatives of \(y_p\), their Hilbert transforms,
and capillary contributions involving \(T_{1,p}\) and \(T_{2,p}\).
Each scale is normalized by its grid average and clipped to prevent
extreme reweighting, while \(S_{\rm der}\) places greater emphasis on
the curvature of the base fields. The scales are recomputed whenever
the grid is changed and are held fixed during parameter
differentiation, preventing the correction network from reducing the
loss merely by altering its normalization.

With component-specific collocation weights \(\omega_{i\ell}\) and
derivative weights \(\nu_j\), the refinement loss is
\begin{equation}
\begin{aligned}
\mathcal L={}&
 w_{\rm eq}\sum_{u_i\in\mathcal C}\sum_{\ell=1}^{3}
 \omega_{i\ell}
 \left|\frac{e_\ell(u_i)}{S_{\rm eq}(u_i)}\right|^2
+w_{\rm der}\sum_{v_j\in\mathcal D}\nu_j
 \left|\frac{\partial_u e_1(v_j)}{S_{\rm der}(v_j)}\right|^2
 +\mathcal L_{\rm aux}.
\end{aligned}
\end{equation}
The derivative term was strengthened
for the final correction in order to suppress rapidly varying
residual errors. The weights \(\omega_{i\ell}\) emphasize the
periodic endpoints and, in intermediate passes, selected residual
peaks; they are normalized so that the total weight assigned to each
residual component remains fixed. The derivative weights \(\nu_j\)
are uniform. The term \(\mathcal L_{\rm aux}\) contains the retained
amplitude and crest--trough separation penalties. The derivative term
acts on \(e_1\), the main equation evaluated with the current
auxiliary fields, and agrees with
\(\partial_u\mathcal B(y_P)\) when their defining residuals vanish.

After the final refinement pass, the corrected fields were evaluated
on shifted grids before Fourier projection and Newton polishing.
Rigorous certification was then carried out separately on the
resulting frozen trigonometric polynomial.

\begin{figure}[htbp]
    \centering
    \includegraphics[width=\linewidth]
    {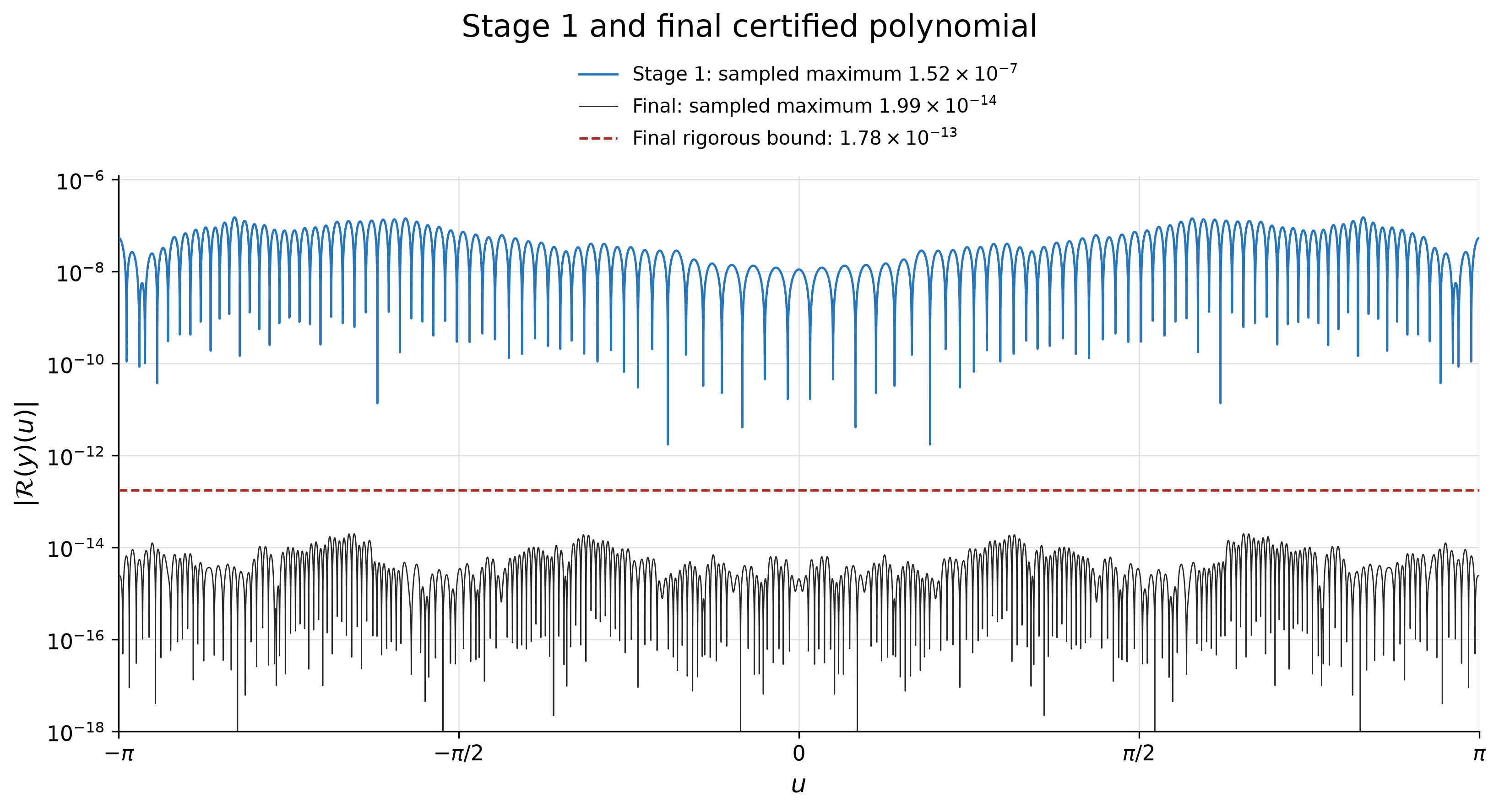}
    \caption{
    Reduction of the scalar Babenko residual.
    Absolute residual as a function of the conformal coordinate
    \(u\) for the first-stage neural profile (blue) and the final
    trigonometric polynomial (black), obtained after successive
    neural corrections, Fourier projection, and Newton polishing.
    Both curves evaluate the scalar Babenko equation with the
    auxiliary fields eliminated, at the physical parameters of
    the corresponding profile. The curves show pointwise numerical
    evaluations; the dashed red line gives the rigorous uniform
    residual bound for the final polynomial.
    }
    \label{fig:residual-refinement}
\end{figure}

Figure~\ref{fig:residual-refinement} compares the scalar Babenko
residual of the first-stage profile with that of the final
trigonometric polynomial, illustrating the reduction achieved through
successive neural corrections, Fourier projection, and Newton
polishing. The rigorous uniform bound for the final polynomial is
shown separately.

\subsection{Computer-assisted accuracy proof}
We certify a uniform bound on the Babenko residual of a fixed finite profile, as follows:

Let 
\begin{equation}
  y_N(u)=\frac{a_0}{2}+\sum_{k=1}^{N}a_k\cos(ku),
  \label{eq:methods-certified-cosine}
\end{equation}
be a fixed finite even cosine series and let $\mu,\beta>0$ be given rational parameters. Define the residual of the Babenko equation with surface tension by
\begin{equation}
\begin{split}
 \mathcal{R}(y)=-\mathcal B(y)={}&\mu\hat H[y_u]-y-y\hat H[y_u]-\hat H[yy_u]
 \\& +\beta\,\partial_u\!\left(\frac{y_u}{D}\right)
 -\beta\hat H\left[
        \partial_u\left(
        \frac{1+\hat H[y_u]}{D}
        \right)
    \right],
\end{split}
\label{eq:supp-cap-residual}
\end{equation}

where
\begin{equation}
    D(u)=
    \sqrt{
        y_u^2+
        \left(1+\hat H[y_u]\right)^2
    },
    \label{eq:supp-D}
\end{equation}
and $\hat H$ denotes the $2\pi$-periodic Hilbert transform, with convention
\begin{equation}
    \hat H[\cos(ku)]=\sin(ku),
    \qquad
    \hat H[\sin(ku)]=-\cos(ku),
    \qquad k\geq1.
    \label{eq:supp-H-convention}
\end{equation}
The goal is to establish $D(u)\geq d_*>0$ for every $u\in\mathbb R$ and $\|\mathcal R(y_N)\|_{L^\infty(\mathbb R)}\leq\epsilon$ for explicit certified constants $d_*$ and $\epsilon$.

\vspace{2mm}
The rigorous proof is split between Julia and Lean 4. A Julia verifier performs the substantive computation using rigorous interval ball arithmetic: it treats the final computed wave as a fixed trigonometric polynomial and establishes a uniform bound on its Babenko residual together with a positive lower bound on the capillary denominator, all with certified error enclosures. A Lean development provides the logical scaffolding: it formalizes the relevant operators and reusable estimates, imports the two final numerical conclusions as explicit assumptions, and checks the deduction of the paper’s residual theorem. The Lean kernel checks this conditional deduction without recomputing the numerical enclosures. Once the polynomial has been frozen, the network weights and training trajectory play no role in the proof. Source code, archived run evidence, build instructions, and the supplementary mathematical account are indexed in the \href{https://github.com/michalshavitNYU/Rigorously-certified-water-waves}{paper repository}.

\subsubsection{Computer-assisted interval certification of the numerical residual}
\label{sec:julia-interval-design}

After neural training and continuation, we project the profile onto
an even cosine series and polish its coefficients with a
finite-dimensional Newton iteration. The resulting polynomial has
the convention in Eq.~\eqref{eq:methods-certified-cosine}. Its stored
binary coefficients and parameters are interpreted as exact dyadic
rationals and fixed before validation. The certified statement
concerns this explicit polynomial and these exact parameter values.
From this point onward, the certification is independent of the
neural network. The normalization and computational settings for
the reported profile are collected in
Section~\ref{sec:certified-profile-results}.

Julia and Arb compute outward-rounded interval enclosures for the
nonlinear residual, add analytic corrections for unresolved Fourier
modes, and use Taylor's theorem to pass from a grid to the whole
torus.

The Julia package \texttt{BabenkoRigorous.jl} encloses every
numerical operation with outward-rounded Arb balls. It first
constructs finite Fourier representations of \(y_N\),
\(y_{N,u}\), \(\hat H[y_{N,u}]\) and their polynomial products.
Differentiation and the periodic Hilbert transform act mode by
mode, while multiplication uses finite coefficient convolution,
with no projection of newly generated modes back to degree \(N\).

Nonlinear capillary quantities require a second representation.
With \(D\) as defined in Eq.~\eqref{eq:supp-D}, the verifier forms
\begin{equation}
  G(u)=\frac{1+\hat H[y_{N,u}](u)}{D(u)}.
  \label{eq:methods-DG}
\end{equation}
At each node of a periodic grid, a Taylor record contains interval
enclosures of a quantity and its derivatives, together with global
derivative suprema. Product, reciprocal and square-root rules
propagate these records. A grid lower scan combined with off-grid
derivative control first proves \(D(u)\geq d_*>0\) over the entire
period; only then are \(D^{-1}\) and the capillary fractions
evaluated.

The Hilbert-transform work is divided into two pathways.
\begin{description}
  \item[Path A: finite trigonometric algebra.]
  The input has a finite Fourier representation, so \(\hat H\)
  acts mode by mode using Eq.~\eqref{eq:supp-H-convention}.
  No sampling or aliasing correction is needed for this operation.

  \item[Path B: certified sampled transform.]
  The input is represented by interval-enclosed values and
  derivatives on a uniform grid, rather than by a finite
  coefficient vector. A discrete Hilbert transform is evaluated
  by certified convolution, and a rigorous Fourier-tail correction
  is added to enclose the continuous transform.
\end{description}
Path B is essential for \(\hat H[\partial_uG]\), because the
square root and reciprocal in Eqs.~\eqref{eq:supp-D} and
\eqref{eq:methods-DG} generally give \(G\) an infinite Fourier
series. The current Julia implementation also routes
\(\hat H[y_Ny_{N,u}]\) through Path B for uniformity, although
\(y_Ny_{N,u}\) is itself a finite trigonometric polynomial of
degree at most \(2N\).

If \(\|F^{(p)}\|_\infty\leq C_p\) for \(p\geq2\), the correction
at each grid node is bounded by
\begin{equation}
  \varepsilon_{\mathrm{alias}}(F)
  =\frac{4C_p}
  {(p-1)(N_{\mathrm{grid}}/2-1)^{p-1}}.
  \label{eq:methods-alias}
\end{equation}
Arb encloses arithmetic roundoff; the aliasing correction is a
separate analytic estimate supplied by the verifier.

The full residual is assembled from these interval records using
the same terms as Eq.~\eqref{eq:supp-cap-residual}. Global derivative
bounds are obtained from the finite Fourier coefficients for
trigonometric polynomials and propagated through products,
reciprocals and square roots using differentiation rules and
interval arithmetic. The positive lower bound on \(D\) controls
the denominators arising in these estimates. For sampled
Hilbert-transform terms, Fourier-decay estimates provide the
additional derivative bounds required by the Taylor remainder.

For grid nodes \(u_j\) separated by \(h=2\pi/N_{\mathrm{grid}}\),
let \(I_{j,m}\) enclose the \(m\)-th derivative of the residual
at \(u_j\), and set \(U_{j,m}=\sup_{x\in I_{j,m}}|x|\).
A global bound
\(C_{M+1}\geq\|\partial_u^{M+1}\mathcal B(y_N)\|_\infty\)
gives, for \(u_j\leq u\leq u_j+h\),
\begin{equation}
  |\mathcal B(y_N)(u)|
  \leq \sum_{m=0}^{M}\frac{U_{j,m}h^m}{m!}
       +\frac{C_{M+1}h^{M+1}}{(M+1)!}.
  \label{eq:methods-cell-bound}
\end{equation}
The Taylor order \(M\) is chosen together with the grid size
\(N_{\mathrm{grid}}\) so that the remainder bound is sufficiently
small to achieve the desired certification accuracy.
Taking the outward-rounded maximum over all cells then proves
a uniform bound over the entire period. The denominator scan
uses the corresponding lower-bound argument.

The production computation uses a sharded implementation: local
Taylor data are streamed to disk, Hilbert-transform derivative
columns are processed separately, and the final residual is
scanned shard by shard.

The implementation, input coefficients, run logs and reproduction
commands are in the
\href{https://github.com/michalshavitNYU/Rigorously-certified-water-waves/tree/main/Julia}
{Julia directory}. The detailed interval estimates are in the
\href{https://github.com/michalshavitNYU/Rigorously-certified-water-waves/blob/main/docs/supplementary_computer_assisted_accuracy_proof.pdf}
{Supplementary Methods}.

\subsubsection{Formal verification in Lean}
\label{sec:lean-certificate-methods}
The Lean formalization provides the mathematical verification layer of a larger computer-assisted proof. Julia/Arblib acts as the certificate generator: it validates a fixed finite cosine-series profile $y_N$, evaluates its Babenko residual $\mathcal{R}(y_N)$ on a fine grid, and uses interval arithmetic and analytic estimates to control rounding, Fourier aliasing and truncation, and between-grid errors. It exports exact rational coefficients, equation parameters and certified endpoints, together with the grid size and Taylor order. The residual endpoint $\varepsilon_0$ bounds the continuous residual, including the Taylor remainder.

For the reported candidate, Lean takes two explicit numerical hypotheses from the validated computation: $D(u)\geq d_*> 0$ for every real $u$, and $|\mathcal{R}(u)|\leq \varepsilon_0$ on $[-\pi,\pi]$. From these inputs, Lean proves $\|\mathcal{R}\|_{L^\infty(\mathbb{R})} \leq \varepsilon_0$ by periodic extension and checks the exact rational comparisons yielding the reported constants. These two hypotheses form the numerical trust boundary: the validated computation supplies the numerical enclosures, while Lean verifies the mathematical deductions. In the Lean implementation, the one-period numerical hypothesis is stated on $[0,2\pi]$; by $2\pi$-periodicity, this is equivalent to the formulation on $[-\pi,\pi]$ used here.

Independently of this numerical instance, a reusable Lean library formalizes the Fourier definition of the Hilbert transform $H$, its action on trigonometric modes and commutation with differentiation, aliasing estimates, product and composition derivative bounds, Taylor-record semantics, and the Taylor grid-to-continuum argument. The publication proof contains no \texttt{sorry} placeholders or \texttt{sorryAx} dependency.

The Lean source code, certificate data and build instructions are
available in the \href{https://github.com/michalshavitNYU/Rigorously-certified-water-waves/tree/main/Lean}
{Lean directory} of the
\href{https://github.com/michalshavitNYU/Rigorously-certified-water-waves/tree/main}
{Paper repository}.
Table~\ref{tab:lean-modules} summarizes the Lean modules and their roles.
Their mathematical content is described in the 
\href{https://github.com/michalshavitNYU/Rigorously-certified-water-waves/blob/main/docs/supplementary_computer_assisted_accuracy_proof.pdf}
{Supplementary Methods}.


\begingroup
\small
\setlength{\tabcolsep}{4pt}
\renewcommand{\arraystretch}{1.13}
\begin{longtable}{@{}p{0.34\textwidth}p{0.48\textwidth}p{0.13\textwidth}@{}}
\caption{Lean modules and their roles in the certification framework.}\label{tab:lean-modules}\\
\toprule
File(s) & Contents & Role \\
\midrule
\endfirsthead
\toprule File(s) & Contents & Role \\ \midrule
\endhead
\bottomrule
\endfoot
\multicolumn{3}{@{}l}{\textbf{Definitions and reusable mathematical library}}\\
\lf{Definitions} & Defines the Fourier-multiplier Hilbert transform, finite profile, derivative, denominator $D$, capillary expression $G$, and residual $R$. & Definitions \\
\lf{HilbertOnTrigPoly} & Proves the Hilbert-transform identities for finite trigonometric profiles, identifying $H[y_u]$ with the explicit cosine sum \texttt{Hyu}. & Analysis \\
\lf{FourierTailHilbert} & Fourier-mode and tail estimates, including the $\pi^2/3$ Hilbert-transform bound. & Analysis \\
\lf{HilbertSmoothness}; \lf{HilbertSharper}; \lf{HilbertCommutesDeriv} & Regularity and derivative commutation for the Hilbert transform under the stated smoothness and periodicity hypotheses. & Analysis \\
\lf{AliasingErrorBound} & Sampled DFT aliasing estimate and equivalence with the cotangent-kernel convolution. & Analysis \\
\lf{LeibnizBound}; \lf{FaaDiBrunoBound} & Derivative bounds for products and compositions. & Analysis \\
\lf{TaylorRemainder} & Taylor-polynomial and remainder estimates on an interval. & Analysis \\
\lf{GridToSupNormTaylor} & Converts cellwise Taylor bounds into a uniform bound on the full grid interval. & Analysis \\
\lf{TaylorSeries\br Correctness} & Defines Taylor records and their semantic invariant; proves arithmetic and grid-to-continuum contracts. & Definitions and proofs \\
\lf{SeriesCorrectness} & Correctness results for square-root and reciprocal derivative records under explicit hypotheses. & Conditional correctness theorems. \\
\lf{TightenSupBoundsLemma} & Justifies downward tightening of derivative bounds. & Proved analysis \\
\lf{FFTHilbert\br Correctness} & Correctness contract for Hilbert-transform Taylor records, conditional on the required enclosures and analytic bounds. & Conditional correctness theorem \\
\lf{DLowerBound} & Deduces positivity from a positive lower bound and records the quotient identity for $G$. & Proved implications \\
\lf{PipelineSmoothness}; \lf{RSmoothness} & Smoothness and periodicity of profile and residual components, with denominator positivity where required. & Analysis \\
\lf{DminBoundCorrectness}; \lf{TightDGSupBounds} & Denominator-bound correctness conditions and derivative bounds for the denominator/capillary chain. & Conditional correctness theorem \\
\lf{ResidualAssembly\br Correctness} & Generic residual assembly; periodicity of $R$ and extension of a one-period bound to all real arguments. & Conditional correctness theorem and analysis \\
\lf{Common}; \lf{Foundations} & Shared foundational imports and support for the library. & Infrastructure \\
\midrule
\multicolumn{3}{@{}l}{\textbf{Concrete publication certificate}}\\
\lf{PaperCertificate} & Exact rational profile coefficients, parameters, grid/Taylor metadata, and numerical certificate endpoints. & Data definitions \\
\lf{PaperCertificateAxioms} & Declares two numerical assumptions: global positive denominator lower bound and residual bound on one period. Proves consequences including periodic extension. & Numerical assumptions and deductions \\
\lf{PaperTheorems} & Publication residual and denominator bounds; exact comparisons with the reported decimal constants. & Proved from certificate assumptions \\
\lf{PaperAudit} & Checks the publication declarations and axiom dependencies; proves the sign-convention bridge to the manuscript residual. & Audit and proofs \\
\midrule
\end{longtable}
\endgroup

\subsubsection{Certified profile, computational settings and bounds}
\label{sec:certified-profile-results}

For the candidate reported in the main text, the certified input is the
polynomial in Eq.~\eqref{eq:methods-certified-cosine}, with stored binary
coefficients interpreted as exact dyadic rationals. In the normalization
$g=1$, the certified parameters are the stored binary64 values:
\begin{equation}
    \mu_0=\frac{1538048895284245}{1125899906842624},
    \qquad
    \beta_0=\frac{2233821443972785}{9007199254740992}.
\end{equation}
The coefficient and parameter values are supplied with the archived input
in the
\href{https://github.com/michalshavitNYU/Rigorously-certified-water-waves}
{paper repository}.
The polynomial has degree $N=150$. The validated calculation uses
\begin{equation}
    N_{\mathrm{grid}}=2{,}000{,}000,
    \qquad
    M=6,
\end{equation}
at $128$-bit Arb precision.

For these fixed inputs, the computation yields
\begin{equation}
    \boxed{
    \|\mathcal{B}_{\mu_0,\beta_0}(y_{150})\|_{L^\infty(\mathbb{R})}
    \leq 1.78\times10^{-13}}
    \label{eq:supp-final-residual}
\end{equation}
and
\begin{equation}
    \boxed{
    D(u)\geq0.3446588>0
    \qquad
    \text{for every }u\in\mathbb{R}.}
    \label{eq:supp-final-D}
\end{equation}
The exact rational certificate endpoints underlying these displayed
bounds are recorded in
\path{Lean/WaveCAP/PaperCertificate.lean}
in the paper repository.

The bound \eqref{eq:supp-final-residual} accounts for propagated Arb
rounding errors, discrete-Hilbert aliasing errors, nonlinear derivative
bounds, and the Taylor remainder controlling the intervals between
grid points. It is therefore a rigorous uniform residual certificate
for the explicit polynomial $y_{150}$ at the stated exact parameter values.

\subsection{Reproducibility and file map}
\label{sec:supp-cap-reproducibility}

The inputs, logs and certificate artifacts for the certified run
are archived in
\begin{center}
\path{Julia/results/phase4_scalar_n150_20260827_2m_20260829/}
\end{center}
within the paper repository. These include the fixed input
coefficients and parameters, reproduction scripts, Arb endpoint
dumps, and the generated Lean certificate, together with checksums.
 The certificate-specific
theorems and assumption audit were checked with the pinned
Lean and Mathlib versions recorded in the repository.

The certified interval calculation was performed on an NYU CIMS
compute server. Table~\ref{tab:certification-runtime} summarizes
the execution times and memory use recorded in the archived logs.
The eight workers ran concurrently, so their wall times are not
summed. The slowest worker followed by the three serial stages
gives an estimated production wall time of approximately
\(2\) hours \(50\) minutes after warmup. The sum of the individual
worker memory peaks was \(61.5\) GB; this does not represent a
measured simultaneous peak for the full run.

\begin{table}[h!]
  \centering
  \footnotesize
  \setlength{\tabcolsep}{4pt}
  \renewcommand{\arraystretch}{1.15}
  \caption{\textbf{Execution time and memory use.}
  Wall times are taken from the system timing logs.
  Peak RSS denotes peak resident memory per process.}
  \label{tab:certification-runtime}
  \begin{tabular}{lrrr}
    \toprule
    Stage & Threads & Wall time & Peak RSS (GB) \\
    \midrule
    Worker 1 & 4 & 1 h 57 min 25 s & 5.019 \\
    Worker 2 & 4 & 2 h 14 min 58 s & 9.242 \\
    Worker 3 & 4 & 1 h 56 min 25 s & 4.972 \\
    Worker 4 & 4 & 2 h 07 min 07 s & 9.287 \\
    Worker 5 & 4 & 2 h 13 min 04 s & 9.324 \\
    Worker 6 & 4 & 1 h 58 min 34 s & 5.088 \\
    Worker 7 & 4 & 2 h 09 min 00 s & 9.267 \\
    Worker 8 & 4 & 2 h 13 min 52 s & 9.283 \\
    \midrule
    Merge & 8 & 1 min 22.08 s & 0.481 \\
    Hilbert-transform stage & 8 & 21 min 33.92 s & 5.149 \\
    Final residual scan (\texttt{s3})
      & 8 & 11 min 49.78 s & 0.477 \\
    \bottomrule
  \end{tabular}
\end{table}

Table~\ref{tab:certification-final-log} reproduces the numerical
summary printed by the final \texttt{s3} stage. The displayed
digits summarize the archived Arb enclosures; the Lean certificate
uses exact rational endpoints. The top-level run log records
\texttt{ALL\_STAGES\_COMPLETE}.

\begin{table}[h]
  \centering
  \footnotesize
  \renewcommand{\arraystretch}{1.2}
  \caption{\textbf{Final-stage numerical summary.}
  Values transcribed from the archived \texttt{s3} log.}
  \label{tab:certification-final-log}
  \begin{tabular}{ll}
    \toprule
    Log entry & Reported value \\
    \midrule
    Bound midpoint
      & \(+1.9928311263803702\times10^{-14}\) \\
    Bound radius
      & \(1.5799937398052128\times10^{-13}\) \\
    Bound interval: lower endpoint
      & \(-1.3807106271671758\times10^{-13}\) \\
    Bound interval: upper endpoint
      & \(+1.7792768524432498\times10^{-13}\) \\
    Residual upper bound
      & \(1.7792768524432498\times10^{-13}\) \\
    Denominator lower bound
      & \(3.4465883518530394\times10^{-1}\) \\
    \bottomrule
  \end{tabular}
\end{table}

\section{Data and code availability}\label{sec:data}\label{sec:code}

The Julia interval-arithmetic code, numerical certificate data, archived
certification run records and Lean formalization are available in the
paper repository, together with pinned software environments and
reproduction instructions:
\url{https://github.com/michalshavitNYU/Rigorously-certified-water-waves}.

\backmatter

\section*{Declarations}

\noindent\textbf{Declaration of generative AI use.} The authors wrote and ran the physics-informed neural-network (PINN) code and wrote the main manuscript. Large language models provided substantial assistance in developing the Julia interval-arithmetic implementation and Lean formalization, and assisted with language editing and proofreading of the manuscript. The numerical certificate relies on Julia/Arb interval arithmetic, while Lean checks the formal deductions conditional on the explicitly stated numerical assumptions.

\noindent\textbf{Competing interests} The authors declare no competing interests.



\providecommand{\CandidateDegree}{150}
\providecommand{\CandidateCount}{151}
\providecommand{\VerificationGrid}{2{,}000{,}000}
\providecommand{\TaylorOrder}{6}
\providecommand{\WorkingPrecision}{128}

\bibliography{sn-bibliography}


\end{document}